\documentclass[a4paper]{article}

\usepackage[numbers,sort&compress]{natbib}
\usepackage{epsfig}
\usepackage{epstopdf}
\usepackage{flushend} 
\usepackage{amsmath}
\usepackage{amsthm}
\usepackage{amssymb}
\usepackage[table]{xcolor}
\usepackage{optidef}
\usepackage{tikz}

\usetikzlibrary{math} 
\usetikzlibrary{shapes}
\usepackage{tkz-euclide}
\usepackage{xstring}

\def\splicelist#1{
\StrCount{#1}{,}[\numofelem]
\ifnum\numofelem>0\relax
     \StrBehind[\numofelem]{#1}{,}[\mylast]%
\else
    \let\mylast#1%
\fi
}

\newcommand{\myroundpoly}[3][thin,color=black]{
\splicelist{#2}
\foreach \vertex [remember=\vertex as \succvertex
    (initially \mylast)] in {#2}{
    \coordinate (\succvertex-next) at ($(\succvertex)!#3!90:(\vertex)$);
    \coordinate (\vertex-previous) at ($(\vertex)!#3!-90:(\succvertex)$);
    \draw[#1] (\succvertex-next) --  (\vertex-previous);
}
\foreach \vertex in {#2}{
    \tkzDrawArc[#1](\vertex,\vertex-next)(\vertex-previous)
}
}
\makeatletter
\def\@normalsize{\@setsize\normalsize{12pt}\xpt\@xpt
\abovedisplayskip 10pt plus2pt minus5pt\belowdisplayskip \abovedisplayskip
\abovedisplayshortskip \z@ plus3pt\belowdisplayshortskip 6pt plus3pt
minus3pt\let\@listi\@listI}

\def\subsize{\@setsize\subsize{12pt}\xipt\@xipt}
\def\section{\@startsection {section}{1}{\z@}{24pt plus 2pt minus 2pt}
{12pt plus 2pt minus 2pt}{\large\bf}}
\def\subsection{\@startsection {subsection}{2}{\z@}{12pt plus 2pt minus 2pt}
{12pt plus 2pt minus 2pt}{\subsize\bf}}
\makeatother

\newtheorem{theorem}{Theorem}

\newtheorem{lemma}{Lemma}
\newtheorem{proposition}{Proposition}
\theoremstyle{remark}
\newtheorem*{remark}{Remark}

\begin{document}

\date{}

\title{\Large\bf A note on the Tuza constant $c_k$ for small $k$
}

\author{Yun-Shan~Lu$^1$, Hung-Lung Wang$^1$}

\date{
$^1$Department of Computer Science and Information Engineering \\
National Taiwan Normal University, Taipei, Taiwan   \\
\{60947097s,hlwang\}@ntnu.edu.tw\\
}

\maketitle


\begin{abstract}
 For a hypergraph $H$, the transversal is a subset of vertices whose intersection with every edge is nonempty. The cardinality of a minimum transversal is the transversal number of $H$, denoted by $\tau(H)$. The Tuza constant $c_k$ is defined as $\sup{\tau(H)/ (m+n)}$, where $H$ ranges over all $k$-uniform hypergraphs, with $m$ and $n$ being the number of edges and vertices, respectively. We give an upper bound and a lower bound on $c_k$. The upper bound  improves the known ones for $k\geq 7$, and the lower bound improves the known ones for  $k\in\{7, 8, 10, 11, 13, 14, 17\}$.    
\end{abstract}

\section{Introduction}

For a hypergraph $H$, the \textit{transversal} is a subset of vertices that intersects every edge. The cardinality of a minimum transversal is the \textit{transversal number} of $H$, denoted by $\tau(H)$. 
The transversal number is a fundamental quantity investigated in hypergraph theory~\cite{HenningYeo2020}. It draws considerable attention since it generalizes several concepts, e.g. the domination number of a graph~\cite{DorflingHenning15,HenningYeo11} and the size of a minimum vertex cover~\cite{ErdosTuza92}. In general, the transversal number is  hard to determining so bounds on the number are sought. The investigation of \textit{Tuza constants} is one direction for this purpose.   A hypergraph $H$ is $k$-uniform if the cardinality of every edge is $k$. Let $\mathcal{H}_k$ be the family of $k$-uniform hypergraphs. Then the \textit{Tuza constant} $c_k$ is defined as 
\[c_k=\sup_{H\in \mathcal{H}_k} \frac{\tau(H)}{m+n},\]
where $m$ and $n$ are the number of edges and vertices of $H$, respectively. { The problem of determining the constant for every $k$ was proposed by Tuza~\cite{TUZA1990117}.} A general upper bound on $c_k$ was proposed by Alon~\cite{Alon90}, which states the following.

\begin{proposition}[Alon~\cite{Alon90}]
\label{prop:ubd_alon}
	Let $H$ be a $k$-uniform hypergraph, where $k>1$, with $n$ vertices and $m$ edges. For any positive real $\alpha$
	\[\tau(H)\leq n\frac{\alpha\ln k}{k}+\frac{m}{k^{\alpha}}.\]
\end{proposition}

\medskip
By Proposition~\ref{prop:ubd_alon} an upper bound on $c_k$ is obtained for $k>1$, namely
\begin{equation}\label{eq:ubd_alon}
    c_k \leq \max \left\{\frac{\alpha\ln k}{k},\frac{1}{k^{\alpha}}\right\},
\end{equation}
where $\alpha$ is any positive real. By taking $\alpha=1$ we have $c_k\leq \ln k/k$. Notice that the minimum is achieved when $\alpha\ln{k}\cdot k^\alpha =k$. Thus, in the following we take $\alpha\ln{k} = W_0(k)$, where $W_0$ is the principle branch of  Lambert's $W$ function~\cite{CGHJK96}. 

\medskip
 For a general lower bound, Lai and Chang~\cite{LaiChang90} provided a way to construct a $k$-uniform hypergraph $H $ for $k\geq 1$, with $\tau(H)=2$ and $m+n=k+1+\lfloor\sqrt{k}\rfloor + \lceil k/\lfloor\sqrt{k}\rfloor\rceil$. Thus, 
\begin{equation}\label{eq:lbd_lai_chang}
    c_k \geq \frac{2}{k+1+\lfloor\sqrt{k}\rfloor + \lceil k/\lfloor\sqrt{k}\rfloor\rceil}.
\end{equation}
We note that the bound of Lai and Chang is exactly $c_k$ for $k\in\{1, 2, 3, 4\}$. Alon~\cite{Alon90} also showed that $c_k \geq (1-O(\ln\ln k/\ln k))\ln k/k$, which implies that the upper bound given in Eq.~\eqref{eq:ubd_alon} is asymptotically tight. However, the exact value for specific $k$ remains unknown, except for $k\leq 4$~\cite{ErdosTuza92,LaiChang90,ChvatalMcDiarmid92,TUZA1990117}. For $k\in\{5,6\}$, the following bounds have been shown~\cite{HenningYeo11,DorflingHenning15,HenningYeo2020}: 
\[\frac{3}{16}\leq c_5\leq \frac{5}{22} \quad\mbox{and}\quad c_6\leq \frac{2569}{14145}.\]

\medskip
In this paper, we improve the upper bounds given in Eq.~\eqref{eq:ubd_alon} for $k\geq 7$, with the idea given by Henning and Yeo~\cite{HenningYeo11}. For the lower bounds, we improve the bounds for $k\in\{7, 8, 10, 11, 13,14, 17\}$ by giving a specific instance. The main result of this paper is given in Theorem~\ref{thm:main}.

\begin{theorem}\label{thm:main}
    For $k\geq 7$ 
    \[c_k \leq \left(1-\frac{1}{2k-1}\right) \frac{W_0(k)}{k}.\]
    In addition, $c_k$ is bounded below by 
    \[\frac{3}{2k+6}.\]
\end{theorem}

\bigskip
A comparison is given in Table~\ref{tab:results}. We prove Theorem~\ref{thm:main} by giving the upper and lower bounds in Sections~\ref{sec:ubd} and~\ref{sec:lbd}, respectively.  

\begin{table}
\centering
\caption{The comparison with known lower bounds (LB) and upper bounds (UB). Marked ones are improvements. }
    \label{tab:results}
    \vspace{0.2cm}
\begin{tabular}{c cccc} 
 \hline
 $k$ & known LB (Eq.~\eqref{eq:lbd_lai_chang})  & the proposed LB  & the proposed UB & known UB (Eq.~\eqref{eq:ubd_alon})  \\ [0.5ex] 
 \hline\hline
 7 & $.1428$ & \cellcolor{blue!10}$.1500$ & \cellcolor{blue!10}$.1762$ & .2178\\ 
 8 & $.1333$ & \cellcolor{blue!10}$.1363$ & \cellcolor{blue!10}$.1634$ & .2008\\
 9 &  $.1250$ & $.1250$ & \cellcolor{blue!10}$.1531$ & .1866\\
 10 & $.1111$ & \cellcolor{blue!10}$.1153$  \cellcolor{blue!10}& \cellcolor{blue!10}$.1438$ & .1746\\
 11 & $.1052$ & \cellcolor{blue!10}$.1071$ & \cellcolor{blue!10}$.1361$ & .1642\\  
 12 & $.1000$ & $.1000$ & \cellcolor{blue!10}$.1292$ & .1553\\
 13 & $.0909$ & \cellcolor{blue!10}$.0937$ & \cellcolor{blue!10}$.1231$ & .1474 \\
 14 & $.0869$ & \cellcolor{blue!10}$.0882$ & \cellcolor{blue!10}$.1177$ & .1403\\
 15 & $.0833$ & $.0833$ & \cellcolor{blue!10}$.1127$ & .1340\\
 16 & $.0800$ & $.0789$ & \cellcolor{blue!10}$.1083$ & .1284\\
 17 & $.0740$ & \cellcolor{blue!10}$.0750$ &  \cellcolor{blue!10}$.1042$ & .1232\\
 18 & $.0714$ & $.0714$ &  \cellcolor{blue!10}$.1005$ & .1185\\
 \hline
\end{tabular}
\end{table}
\section{Upper bounds on $c_k$}
\label{sec:ubd}

Let $H=(V,E)$ be a hypergraph. The \textit{degree} of a vertex $v$ is the number of edges that contain $v$, denoted by $\deg(v)$. The \textit{maximum degree} of $H$, denoted by $\Delta(H)$, is defined as $\max\{\deg(v)\colon\, v\in V\}$. If every vertex of $H$ has degree $r$, then $H$ is said to be \textit{$r$-regular}. Two vertices are \textit{neighbors} with each other if they belong to the same edge, and the \textit{neighborhood} of a vertex $u$ is the set of all neighbors of $u$, denoted by $N(u)$. Two edges \textit{overlap} if their intersection contains two or more vertices. If $H$ has no edges overlapping, then $H$ is said to be \textit{linear}.

\medskip
Assume that $H$ is a $k$-uniform hypergraph, with $n$ vertices and $m$ edges.  
Let $n_i$ be the number of vertices of degree $i$, and let $n_{\geq i} = \sum_{j\geq i} n_j$.  For a given integer $d$, associate each vertex of degree at least $d$ with a weight $w_d$, each vertex with degree $i$ for $0\leq i<d$ with a weight $w_i$, and each edge with a weight $w_m$. For $k> d \geq 1$, we show that any solution of LP~\eqref{lp:main} gives an upper bound on $c_k$\footnote{We use LP as an abbreviation for Linear Program.}. The idea is extended from Henning and Yeo~\cite{HenningYeo11}. 

\begin{mini!} 
{}{w_m}{\label{lp:main}}{}
\addConstraint{kw_1+ w_m}{\geq 1 \label{cst:1}}
\addConstraint{w_d+(d+1)w_m }{\geq 1 \label{cst:2}}
\addConstraint{ i(k-1)(w_i-w_{i-1})+2w_{i-1}-w_{i-2}+ iw_m}{\geq 1 \label{cst:3}}{\mkern20mu\quad i=2 ,\ldots, d}
\addConstraint{2i(k-1)(w_i-w_{i-1})+w_i+2w_{i-1}-w_{i-2}+ 2iw_m}{\geq 2 \label{cst:4}}{\mkern20mu\quad i=2,\ldots, d}
\addConstraint{w_d-w_{d-1}\leq \dots\leq w_2- w_1}{\leq w_1 \label{cst:5}}{}
\addConstraint{0=w_0\leq w_1\leq \dots\leq w_d}{\leq w_m \label{cst:6}}{}
\end{mini!}


\begin{lemma}\label{lem:upd}
    For $k>d\geq 2$, every $k$-uniform hypergraph $H$ satisfies
    \begin{equation*}
        \tau(H) \leq w(H)\leq w_m(n+m), 
    \end{equation*}
    where $w(H) = \sum_{i=0}^{d-1} w_in_i+n_{\geq d}w_d + mw_m$ and $(w_0,\dots,w_d,w_m)$ is a feasible solution of LP~\eqref{lp:main}.
\end{lemma}

\begin{proof}
    We prove the lemma by induction on $n$. {We may assume that no two edges are identical}.  For $n \leq k$, there is at most one edge. Either $\tau(H)=0$, or by Constraint~\eqref{cst:1}
\[\tau(H) \leq 1 \leq kw_1+w_m= w(H).\] 
 
For $n> k$, we develop the following cases: (i)~$\Delta(H)=1$, (ii) $\Delta(H)\geq d+1$, and (iii) $2\leq \Delta(H)\leq d$. For any subset $U\subseteq V$, the hypergraph $H-U$ is obtained by removing the vertices in $U$ as well as the edges in $\{X\in E\colon\, X\cap U\neq\emptyset\}$. { Notice that a $k$-uniform hypergraph may contain isolated vertices}. Clearly, for any subset $U$ of vertices the hypergraph $H-U$ remains $k$-uniform. When $U$ is a singleton, say $U=\{u\}$, we write $H-U$ as $H-u$ for succinctness. 
    
    \begin{itemize}
    
         \item[-] Case (i): Since $\Delta(H)=1$, the hypergraph consists of isolated vertices and disjoint edges. Let $v$ be a vertex of degree $1$. It follows that
        \[\tau(H) = \tau(H-v) + 1\underset{I.H.}{\leq} w(H-v)+1\underset{\eqref{cst:1}}{\leq} w(H-v) + kw_1 + w_m = w(H).\]
      
        \item[-] Case (ii): Let $v$ be a vertex of $\deg(v)\geq d+1$. Then 
        \begin{equation*}
        w(H-v)\leq w(H) - w_d - (d+1)w_m.
        \end{equation*}
        Along with the inductive hypothesis and Constraint~\eqref{cst:2} we have
         \[\tau(H) \leq  \tau(H-v) +1 \underset{{I.H.}}{\leq} w(H-v)+1 \underset{\eqref{cst:2}}{\leq} w(H-v) + w_d + (d+1)w_m  \leq w(H).\] 
       
        \item[-] Case (iii): Let $d'=\Delta(H)$ with $2\leq d'\leq d$. If there is a vertex $v$ of $\deg(v)=d'$ such that $v$ has a neighbor of degree at most $d'-1$ or there are two edges containing $v$ overlap, then 
       \begin{multline*}
        w(H-v)\leq  w(H) - w_{d'} - d'\cdot w_m \\ - (d' (k-1) (w_{d'}-w_{d'-1}) - (w_{d'}-w_{d'-1})+(w_{d'-1}-w_{d'-2})). 
        \end{multline*} 
        {The inequality holds since with Cosntraint~\eqref{cst:5} there is a neighbor of $v$ having its degree decreased by at least $w_{d'-1}-w_{d'-2}$, and every other neighbor has its degree  decreased by at least $w_{d'}-w_{d'-1}$}. 
        Thus,
        \begin{multline*}
        \tau(H) \leq \tau(H-v) +1 \underset{{I.H.}}{\leq} w(H-v)+1\\ \underset{\eqref{cst:3}}{\leq} w(H-v) + d' (k-1)w_{d'} - (d' (k-1)-2) w_{d'-1} - w_{d'-2}+ d'w_m   \leq w(H). 
        \end{multline*}
        Otherwise, $H$ is $d'$-regular and is linear. Let $v$ be a vertex of $H$. For any neighbor $u$ of $v$, there has to be an edge $e_u$ with $e_u\in E_u\setminus E_v$, where $E_u$ and $E_v$ are the set of edges containing $u$ and $v$, respectively. Since $k> d\geq d'$, there exists $u'\in e_u\setminus N(v)$. Consider the hypergraph $H' = H-\{v, u'\}$, we have 
        \begin{multline*}
        	w(H')\leq w(H) - 2w_{d'} - 2d'w_m \\ - (2d' (k-1) (w_{d'}-w_{d'-1}) - (w_{d'}-w_{d'-1}) + (w_{d'-1}-w_{d'-2})).
        \end{multline*} 
        It follows that
        \begin{multline*}
        \tau(H) \leq  \tau(H') +2\underset{{I.H.}}{\leq} w(H')+2\\ \underset{\eqref{cst:4}}{\leq} w(H') + (2d' (k-1)+1)w_{d'}- (2d' (k-1)-2)w_{d'-1} -w_{d'-2}+2d'w_m \leq w(H). 
        \end{multline*}

    \end{itemize}
    With the analysis above, by mathematical induction the lemma is proved. 
\end{proof}

\begin{remark}
  
Notice that any feasible solution of  LP~\eqref{lp:main} satisfies $(d+2)w_m \geq 1$, relaxed from~\eqref{cst:2}, so the increasing of $d$ provides an improvement of $w_m$ until $(d+2)w_m$ is strictly greater than $1$. To choose a proper value for $d$, one may choose $d=k-1$ to obtain the best possible $w_m$. For $7\leq k\leq 18$, the smallest values for $d$ to achieve the best possible $w_m$ are given in Table~\ref{tab:para}. We note that LP~\eqref{lp:main} can be applied to obtain upper bounds on $c_k$ for $3\leq k\leq 6$, but the resulting bounds are no better than the existing ones.

\end{remark}

\begin{table}[h]
    \centering
    \caption{Parameters and results of LP~\eqref{lp:main}. The upper bounds (UB) stated in Theorem~\ref{thm:main} are rounded up to four decimal places.}
    \label{tab:para}
    \vspace{0.2cm}
    \begin{tabular}{cccc}
        \hline
        $k$ & $d$ & $w_m$ & UB in Theorem~\ref{thm:main}   \\ [1ex]
         \hline\hline 
         $7$ & $4$ & $.1762$ & $.2010$  \\    
         $8$ & $5$ & $.1634$ & $.1874$  \\
         $9$ & $5$ & $.1531$ & $.1756$  \\
         $10$ & $5$ & $.1438$ & $.1654$ \\
         $11$ & $6$ & $.1361$ & $.1563$  \\
         $12$ & $6$ & $.1292$ & $.1485$  \\
         $13$ & $7$ & $.1231$ & $.1415$  \\
         $14$ & $7$ & $.1177$ & $.1351$  \\
         $15$ & $7$ & $.1127$ & $.1293$  \\
         $16$ & $8$ & $.1083$ & $.1242$  \\
         $17$ & $8$ & $.1042$ & $.1194$  \\
         $18$ & $8$ & $.1005$ & $.1151$  \\
         \hline
    \end{tabular}
\end{table}

\subsection{The upper bound stated in Theorem~\ref{thm:main}}
\label{sec:lp_ubd}

In this section we analyze LP~\eqref{lp:main}, showing that the upper bound stated in Theorem~\ref{thm:main} holds. 
 For each $k$ with $k\leq 18$, the bound can be verified by solving LP~\eqref{lp:main}\footnote{ We use SoPlex to compute the exact rational solutions~\cite{GSW16,GS20}.}, as shown in Table~\ref{tab:para}.
In the following, we assume $k\geq 19$ and write $W_0$ as $W$ for succinctness. Some known properties of $W$ are summarized below. First, by definition

\begin{equation}
\label{eq:w_def}
	W(x)e^{W(x)}=x.
\end{equation}

\noindent An upper bound on the principal branch is derived by Hoorfar and Hassani~\cite{Hoorfar2008}: 
For $x\geq -1/e$ and $y>1/e$, 
\begin{equation}\label{eq:ubd_w}
	W(x)\leq \ln\left(\frac{x+y}{1+\ln y}\right),
\end{equation}
where equality holds only if $y\ln y=x$. 

\medskip
Now we are ready to prove the upper bound stated in Theorem~\ref{thm:main} for $k\geq 19$. Observe that Constraints~\eqref{cst:3} and~\eqref{cst:4} can be restricted to be Constraint~\eqref{cst:new} since  $w_i-w_{i-1}\leq w_{i-1}-w_{i-2}$.

\begin{mini!}
{}{w_m}{\label{lp:restricted}}{}
\addConstraint{kw_1+ w_m}{\geq 1 \label{cst:new1}} 
\addConstraint{w_d+(d+1)w_m }{\geq 1 \label{cst:new2}}
\addConstraint{i(k-1)(w_i-w_{i-1})+w_{i}+iw_m}{\geq 1\label{cst:new}}{\mkern20mu\quad i=2 ,\ldots, d} 
\addConstraint{w_d-w_{d-1}\leq \dots\leq w_2- w_1}{\leq w_1 \label{cst:new4}}
\addConstraint{0\leq w_1\leq \dots\leq w_d}{\leq w_m \label{cst:new5}}
\end{mini!}

Let
\begin{equation}\label{eq:wm}
	w_m =  \frac{W(k)}{k} - \epsilon
\end{equation}
with 
\[0<\epsilon\leq \frac{1}{2k-1}\frac{W(k)}{k}.\]
%
%

%
\noindent 
Consider the sequence $(g_0, g_1, \dots)$, where 
\[
	g_i = \begin{cases}
		-1, & {i=0}\\
		\frac{i(k-1)}{i(k-1)+1}g_{i-1}, & i\geq 1.
	\end{cases}
\]
We choose $d$ to be 
\[\max\{x\in\mathbb{N}\colon\, g_x-g_{x-1} \geq {w_m}/{k}\}.\]

\medskip
\noindent
We claim that $(w_1,\dots,w_d,w_m)$ is a feasible solution of LP~\eqref{lp:restricted}, where
\begin{equation}\label{eq:wi}
	w_i = g_i+1-\frac{i}{k} w_m, 
\end{equation}
for $1\leq i\leq d$.  Notice that $d$ is well-defined since 
\[g_1-g_0 = \frac{1}{k}\geq \frac{W(k)}{k^2}\]
for every $k\in\mathbb{N}$.  Moreover, since $w_d-w_{d-1}\geq 0$ and $(g_0, g_1, \dots)$ is an increasing sequence, we have for $d\geq 2$ 
\begin{equation}
\label{eq:d_ubd}
	d<\frac{k}{k-1}\frac{-g_d}{w_m}\leq \frac{k}{k-1}\frac{-g_2}{w_m}=\frac{k}{W(k)}=e^{W(k)}.
\end{equation}
Also notice that for $k\geq 19$, we have $d\geq 4$ since
\[g_4-g_3 = \frac{1}{4k-3}\frac{3k-3}{3k-2}\frac{2k-2}{2k-1}\frac{k-1}{k}\geq \frac{W(k)}{k^2}.\]

\medskip
Clearly, Constraints~\eqref{cst:new1} and~\eqref{cst:new} hold with equality. We claim that Constraints~\eqref{cst:new2},~\eqref{cst:new4} and~\eqref{cst:new5} are all satisfied.

\begin{lemma}
\label{lem:cst6f}
$0\leq w_1\leq \dots\leq w_d$.
\end{lemma}

\begin{proof}
	Clearly, $w_1\geq 0$. 
	For $i\geq 2$, 
	\[w_i-w_{i-1} = g_i-g_{i-1} - \frac{w_m}{k}.\]
	Since $d = \max\{x\in\mathbb{N}\colon\, g_x-g_{x-1} \geq {w_m}/{k}\}$, we have $w_{i}\geq w_{i-1}$ for $2\leq i\leq d$. 
	
\end{proof}

\begin{lemma}
\label{lem:cst6e}
$w_d-w_{d-1}\leq \dots\leq w_2- w_1\leq w_1$.
\end{lemma}

\begin{proof}
	The last inequality holds since
	\[w_2-w_1 = g_2-g_{1} - \frac{w_m}{k} = \frac{k-1}{2k-1}\cdot \frac{1}{k}- \frac{w_m}{k} < \frac{1}{k}- \frac{w_m}{k} = w_1.\]
	For $3\leq i\leq d$
	\begin{multline*}
		w_i-w_{i-1} = g_i-g_{i-1} - \frac{w_m}{k} 
		 		= \frac{-g_{i-1}}{i(k-1)+1}  - \frac{w_m}{k} \\
				\leq  \frac{-g_{i-1}}{(i-1)(k-1)}  - \frac{w_m}{k} 
				 =  g_{i-1}-g_{i-2} - \frac{w_m}{k}
				= w_{i-1}-w_{i-2}.
	\end{multline*}	
\end{proof}

\bigskip
\begin{lemma}
\label{lem:cst6c}
${w_d+(d+1)w_m }{\geq 1}$.
\end{lemma}

\begin{proof}
	Since 
	\[0 < g_d - g_{d+1} +\frac{w_m}{k} = \frac{g_d}{(d+1)(k-1)+1} + \frac{w_m}{k}, \]
	it follows that
	\begin{multline*}
		w_d + (d+1) w_m = g_d + 1 - \frac{d}{k}w_m + (d+1)w_m 
		 >  g_d + 1 - \left(d+1-\frac{d}{k}\right)\frac{kg_d}{(d+1)(k-1)+1}
		= 1.
	\end{multline*}
\end{proof}

\bigskip
Let $H_n$ be the $n$th harmonic number. It is known that for $n\in\mathbb{N}$
\[H_n \leq \ln n + \gamma + \frac{1}{2n},\]
where $\gamma=0.577215\dots$ is the Euler's constant. 

\begin{proposition}
\label{prop:kd_10}
	Let $\mu = \gamma + \frac18$. For $k\geq 19$, 
	\[\frac{2k-2}{2k-1} \cdot \left(1+\ln d\right)
	\geq \frac{k}{k-1} \cdot\left(\mu+\ln d\right).\]
\end{proposition}

\begin{proof}
	Consider the function $f\colon\mathbb{R}\to\mathbb{R}$ with
	\[f(x) = \frac{(x-1)(2x-2)}{x(2x-1)}-\frac{3x -2}{x(2x-1)} \cdot \ln x.\]
	For $x\geq 19$, the function is increasing since
	\[\frac{d}{dx}f(x) = \frac{-x+(6x^2-8x+2)\cdot \ln x}{(2x^2-x)^2} > 0.\]
	In addition, $f(19)\approx 0.6914 > \mu$. Thus, by Eq.~\eqref{eq:d_ubd} and taking $y=k$ in Eq.~\eqref{eq:ubd_w}, we have
	\[\frac{(k-1)(2k-2)}{k(2k-1)}-\frac{3k -2}{k(2k-1)} \cdot \ln d = \frac{(k-1)(2k-2)}{k(2k-1)}-\frac{3k -2}{k(2k-1)} \cdot W(k) \geq f(k) >\mu.\]
\end{proof}

\bigskip
\begin{lemma}
\label{lem:cst6f_right}
For $k\geq 19$,  $w_d\leq w_m$.
\end{lemma}

\begin{proof}
	We show that 
	\[-g_d\geq 1-\left(1+\frac{d}{k}\right)w_m.\]
	

Taking the logarithm of both sides, we have
	\begin{equation*}
		\ln(-g_d) =\ln\prod_{x=1}^d \frac{x}{x+\frac{1}{k-1}} 
		= \ln\prod_{x=1}^d \left(1+\frac{1}{x(k-1)}\right)^{-1} 
		\geq  -\frac{1}{k-1}\sum_{x=1}^d \frac1x
		= -\frac{H_d}{k-1}
	\end{equation*}
	
	\noindent
	and
%
	\begin{equation*}
		\ln \left(1-\left(1+\frac{d}{k}\right)\left(\frac{W(k)}{k}-\epsilon\right)\right) \leq -\left(1+\frac{d}{k}\right)\left(\frac{W(k)}{k}-\epsilon\right),		
	\end{equation*}
	
	\noindent 
	where the inequalities holds due to $e^x \geq 1+x$ for $x\in\mathbb{R}$.  Therefore, it suffices to claim  
	\[\frac{W(k)}{k}-\epsilon \geq \frac{k}{k-1}\cdot\frac{H_d}{k+d}.\]

\noindent 	
Taking $y$ such that $y\ln y=k$ in Eq.~\eqref{eq:ubd_w}, we have
\begin{multline*}
	\frac{W(k)}{k}-\epsilon \geq \frac{2k-2}{2k-1} \cdot \frac{W(k)}{k} = \frac{2k-2}{2k-1} \cdot e^{-W(k)}
	 \underset{\eqref{eq:ubd_w}}{=} \frac{2k-2}{2k-1} \cdot\frac{1+\ln y}{k+y}\\
	\geq  \frac{2k-2}{2k-1} \cdot\frac{1+\ln d}{k+d} 
	\geq \frac{k}{k-1}\cdot\frac{\mu+\ln d}{k+d}\geq \frac{k}{k-1}\cdot\frac{H_d}{k+d},
\end{multline*}

\medskip
\noindent
where the second last inequality follows from Proposition~\ref{prop:kd_10}, and the last one holds since $d\geq 4$ for $k\geq 19$, implying 
\[H_d \leq \gamma + \frac18 +\ln d.\]
Thus, the claim holds, and the lemma is proved.
%
%
%
\end{proof}

\medskip
By Lemmas~\ref{lem:cst6f},~\ref{lem:cst6e},~\ref{lem:cst6c}, and~\ref{lem:cst6f_right}, the tuple $(w_1,\dots,w_d,w_m)$ defined by Eq.~\eqref{eq:wm} and~\eqref{eq:wi} is a feasible solution of LP~\eqref{lp:restricted}. The objective function value of any feasible solution of LP~\eqref{lp:restricted} gives an upper bound on the optimal objective function value of LP~\eqref{lp:main} since the former linear program is  restricted from the latter. Thus, the upper bound stated in Theorem~\ref{thm:main} is proved. 

\section{Lower bounds on $c_k$}
\label{sec:lbd}

\begin{figure}[!t]
\centering
\tikzmath{
    \l1 = 4;
    \l2 = 10;
    \xa = \l1/1.7; 
    \xb = \l2 - (\l2-\l1)/2; 
    \xc = \l2/2;
    \xd = -\xc;
    \xe = -\xb;
    \xf = -\xa;
    \ya = \l1*sqrt(3)/2 + (\l2-\l1)/sqrt(3); 
    \yb = -(\l2-\l1)/(0.75*(2*sqrt(3)));
    \yc = (\l2-\l1)/(2*sqrt(3)) - \l2*sqrt(3)/2;
    \yd = \yc;
    \ye = \yb;
    \yf = \ya;
    \xab = (\xa+\xb)/2;
    \xbc = 0.6*\xb+0.4*\xc;
    \xcd = (\xc+\xd)/2;
    \xde = 0.6*\xd+0.4*\xe;
    \xef = (\xe+\xf)/2;
    \xfa = 0.57*\xf+0.43*\xa;
    \yab = (\ya+\yb)/2;
    \ybc = 0.6*\yb+0.4*\yc;
    \ycd = (\yc+\yd)/2;
    \yde = 0.6*\yd+0.4*\ye;
    \yef = (\ye+\yf)/2;
    \yfa = 0.6*\yf+0.4*\ya;   
    \xz = 0;
    \yz = -0.7;
 } 

\begin{tikzpicture}[scale=0.4]
    \node[circle] (a) at (\xa,\ya){$X_{1,R}$};
    \node[circle] (b) at (\xb,\yb){$X_{2,L}$};
    \node[circle] (c) at (\xc,\yc){$X_{2,R}$};
    \node[circle] (d) at (\xd,\yd){$X_{3,L}$};
    \node[circle] (e) at (\xe,\ye){$X_{3,R}$};
    \node[circle] (f) at (\xf,\yf){$X_{1,L}$};
    \node[circle] (ab) at (\xab,\yab){$Y_{1,2}$};
    \node[circle] (bc) at (\xbc,\ybc){$X_{2,M}$};
    \node[circle] (cd) at (\xcd,\ycd){$Y_{2,3}$};
    \node[circle] (de) at (\xde,\yde){$X_{3,M}$};
    \node[circle] (ef) at (\xef,\yef){$Y_{3,1}$};
    \node[circle] (fa) at (\xfa,\yfa){$X_{1,M}$};
    \node[circle] (z) at (\xz,\yz){$Z$};
    \node[circle] (xx) at (12, 0){};  
    
    \myroundpoly[color = red, ultra thick]{a,b,c,f}{1cm}
    \myroundpoly[color = blue, ultra thick]{a,d,e,f}{1cm}
    \myroundpoly[color = black, dashed, very thick]{b,c,d,e}{1cm}
\end{tikzpicture}
\begin{tikzpicture}[scale=0.4]
    \node[circle] (a) at (\xa,\ya){$X_{1,R}$};
    \node[circle] (b) at (\xb,\yb){$X_{2,L}$};
    \node[circle] (c) at (\xc,\yc){$X_{2,R}$};
    \node[circle] (d) at (\xd,\yd){$X_{3,L}$};
    \node[circle] (e) at (\xe,\ye){$X_{3,R}$};
    \node[circle] (f) at (\xf,\yf){$X_{1,L}$};
    \node[circle] (ab) at (\xab,\yab){$Y_{1,2}$};
    \node[circle] (bc) at (\xbc,\ybc){$X_{2,M}$};
    \node[circle] (cd) at (\xcd,\ycd){$Y_{2,3}$};
    \node[circle] (de) at (\xde,\yde){$X_{3,M}$};
    \node[circle] (ef) at (\xef,\yef){$Y_{3,1}$};
    \node[circle] (fa) at (\xfa,\yfa){$X_{1,M}$};
    \node[circle] (z) at (\xz,\yz){$Z$};
    
    \myroundpoly[color = red, ultra thick]{ab,b,bc,e,ef}{1.1cm}
    \myroundpoly[color = blue, ultra thick]{cd,d,de,a,ab}{1.2cm}
    \myroundpoly[color = black, dashed, very thick]{ef,f,fa,c,cd}{1.1cm}
\end{tikzpicture}
\caption{A $k$-uniform hypergraph with $2k$ vertices and $6$ edges. The transversal number is $3$.   \label{fig:lbd}}
\end{figure}
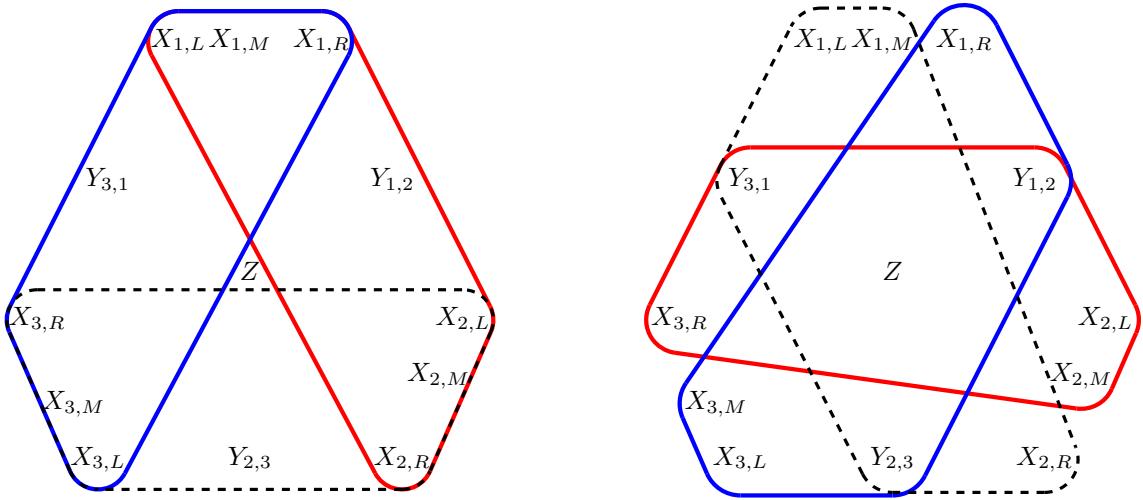

We give the lower bound by defining a $k$-uniform hypergraph $H=(V,E)$ as follows. See Fig.~\ref{fig:lbd} for an illustration. For the notation $[s]$, it denotes the set $\{x\in \mathbb{N}\colon\, 1\leq x\leq s \}$ if $s\in \mathbb{N}$, and stands for the Iverson bracket if $s$ is a logical statement. Let $V$ be the disjoint union of the $13$ sets, $X_{1,L}$, $X_{1,M}$, $X_{1,R}$, $X_{2,L}$, $X_{2,M}$, $X_{2,R}$, $X_{3,L}$, $X_{3,M}$, $X_{3,R}$, $Y_{1,2}$, $Y_{2,3}$, $Y_{3,1}$, and $Z$, whose cardinalities satisfy:

\begin{itemize}
    \item For $i\in[3]$, $|X_{i,L}|=|X_{i,R}|=\left\lfloor\frac{k}{4}\right\rfloor$.
    \item For $i\in[3]$, $|X_{i,M}|=[(k\equiv 2\bmod 4) \vee (k\equiv 3\bmod 4)]$.
    \item $|Y_{1,2}| = |Y_{2,3}| = |Y_{3,1}| = [k\equiv 1\bmod 2]$
    \item $|Z| = k - \left(2\cdot \left\lfloor\frac{k}{4}\right\rfloor + 2\cdot [k\equiv 1\bmod 2] + [(k\equiv 2\bmod 4) \vee (k\equiv 3\bmod 4)] \right)$
\end{itemize}

\noindent
Notice that some of these vertex subsets may be empty, depending on the value of $k$. There are six edges:
\begin{itemize}
    \item $X_{1,L}\cup X_{1,M}\cup X_{1,R}\cup Y_{1,2}\cup X_{2,L}\cup X_{2,M}\cup X_{2,R}$
    \item $X_{2,L}\cup X_{2,M}\cup X_{2,R}\cup Y_{2,3}\cup X_{3,L}\cup X_{3,M}\cup X_{3,R}$
    \item $X_{3,L}\cup X_{3,M}\cup X_{3,R}\cup Y_{3,1}\cup X_{1,L}\cup X_{1,M}\cup X_{1,R}$
    \item $Y_{1,2}\cup X_{2,L}\cup X_{2,M}\cup Z\cup Y_{3,1}\cup X_{3,R}$
    \item $Y_{2,3}\cup X_{3,L}\cup X_{3,M}\cup Z\cup Y_{1,2}\cup X_{1,R}$
    \item $Y_{3,1}\cup X_{1,L}\cup X_{1,M}\cup Z\cup Y_{2,3}\cup X_{2,R}$
\end{itemize}

Any transversal of the hypergraph contains at least two vertices, say $a$ and $b$, such that $\{a, b\}$ intersects all three edges on the left of Fig.~\ref{fig:lbd}. Moreover, such $a$ and $b$ reside in at most two edges on the right of Fig.~\ref{fig:lbd}. Thus,  the transversal number of $H$ is $3$. Since each edge has cardinality $k$, it follows that
\[ \frac{\tau(H)}{|V|+|E|}=\frac{3}{2k+6}.\]

\noindent
As shown in Table~\ref{tab:results}, the ratio is greater than that given in Eq.~\eqref{eq:lbd_lai_chang} for $k\in\{7, 8, 10, 11, 13, 14, 17\}$, which provides improved lower bounds. We note that for $k=5$ the constructed hypergraph is exactly the instance given by Henning and Yeo~\cite{HenningYeo2020}, showing $c_5\geq 3/16$.

\subsection*{Acknowledgements}
The authors would like to thank the anonymous reviewers for helpful comments. Part of the main results were worked out during the revision, based on the reviewers' comments.  This research was supported in part by National Science and Technology Council of Taiwan under contract MOST grant 110-2221-E-003-003-MY3.

\bibliography{ref}

\end{document}